\documentclass[12pt]{amsart}

\usepackage{amsmath}
\usepackage{amsthm}
\usepackage{amssymb}
\usepackage{hyperref}

 \DeclareMathOperator{\Tr}{Tr}

 \DeclareMathOperator{\dimm}{dim}

\usepackage{eucal}
\usepackage{times}
\usepackage{euler}

\begin{document}

 
\title[Permutations without long decreasing subsequences]
{Permutations without long decreasing subsequences and random matrices}
\author{Piotr \'Sniady}
 
\address{Institute of Mathematics,
University of Wroclaw, \mbox{pl.\ Grunwaldzki~2/4,} \mbox{50-384
Wroclaw,} Poland} \email{Piotr.Sniady@math.uni.wroc.pl}
\subjclass[2000]{05E10, 15A52, 60J65}




\begin{abstract}
We study the shape of the Young diagram $\lambda$ associated via the
Robinson--Schensted--Knuth algorithm to a random permutation in
$S_n$ such that the length of the longest decreasing subsequence is
not bigger than a fixed number $d$; in other words we study the restriction of
the Plancherel measure to Young diagrams with at most $d$ rows. We prove that in
the limit $n\to\infty$ the rows of $\lambda$ behave like the eigenvalues of a
certain random matrix (namely the traceless Gaussian Unitary Ensemble random matrix)
with $d$
rows and columns. In particular, the length of the longest
increasing subsequence of such a random permutation behaves
asymptotically like the largest eigenvalue of the corresponding random
matrix.
\end{abstract}

\maketitle

\theoremstyle{plain}
\newtheorem{lemma}{Lemma}
\newtheorem{theorem}[lemma]{Theorem}
\newtheorem{proposition}[lemma]{Proposition}
\newtheorem{corollary}[lemma]{Corollary}
\newtheorem*{conjecture}{Conjecture}
\newtheorem{hypothesis}[lemma]{Hypothesis}
\newtheorem{fact}[lemma]{Fact}
\newtheorem{informal}[lemma]{Informal theorem}
\newtheorem{problem}[lemma]{Problem}

\theoremstyle{definition}
\newtheorem{definition}[lemma]{Definition}

\theoremstyle{remark}
\newtheorem*{remark}{Remark}
\newtheorem*{example}{Example}
\newtheorem{question}[lemma]{Question}

\newcommand{\A}{{\mathfrak{A}}}
\newcommand{\B}{{\mathfrak{B}}}
\newcommand{\E}{{\mathbb{E}}}
\newcommand{\C}{{\mathbb{C}}}
\newcommand{\R}{{\mathbb{R}}}
\newcommand{\Z}{{\mathbb{Z}}}
\newcommand{\El}{{\mathcal{L}}}
\newcommand{\N}{{\mathbb{N}}}
\newcommand{\gwia}{^{\star}}

\newcommand{\M}[1]{M_{#1}(\mathbb{C})}

\newcommand{\RR}{R}

\section{Introduction}

\subsection{Formulation of the problem}
\label{sub:Formula-problem} Let an integer $d\geq 1$ be fixed. For any
integer $n\geq 1$ we consider the set of the permutations $\pi\in
S_n$ such that the length of the longest decreasing subsequence of
$\pi$ is not bigger than $d$; in other words it is the set of the
permutations avoiding the pattern $(d+1,d,\dots,3,2,1)$. Let $\pi_n$
be a random element of this set (probabilities of all elements are
equal). In this article we are interested in the following problem:

\begin{problem}
\label{problem:sigma} Let\/ $\pi_n\in S_n$ be a random permutation
with the longest decreasing subsequence of length at most $d$. What can
we say about the asymptotic behavior of the length of the longest
increasing subsequence of\/ $\pi_n$ in the limit $n\to\infty$?
\end{problem}

Let $\lambda_n=(\lambda_{n,1},\dots,\lambda_{n,d})$ be the (random)
Young diagram associated via the Robinson--Schensted--Knuth
algorithm to $\pi_n$ (notice that since the number of the rows of
$\lambda_n$ is equal to the length of the longest decreasing
subsequence of $\pi_n$, $\lambda_n$ has at most $d$ rows). In
other words, $\lambda_n$ is a random Young diagram with at most $d$
rows, where the probability of the Young diagram $\lambda$ is
proportional to $(\dimm \rho_{\lambda})^2$, where $\dimm
\rho_{\lambda}$ denotes the dimension of the corresponding
irreducible representation of $S_n$; therefore, if we drop the
restriction on the number of the rows of the Young diagrams (which can be
alternatively stated as $d\geq n$), then the distribution of $\lambda_n$ is
the celebrated Plancherel measure.


Since $\lambda_{n,1}$ is equal to the length of the longest
increasing subsequence in $\pi_n$, Problem \ref{problem:sigma}
is a special case of the following more general one:
\begin{problem}
What can we say about the asymptotic behavior of the random
variables $(\lambda_{n,1},\dots,\lambda_{n,d})$ in the limit
$n\to\infty$?
\end{problem}

\subsection{Case $d=2$} The first non-trivial case $d=2$ was considered by
Deutsch, Hildebrand and Wilf \cite{DeutschHildebrandWilf}. In this
case the random variables $\lambda_{n,1},\lambda_{n,2}$ are subject
to a constraint $\lambda_{n,1}+\lambda_{n,2}=n$ therefore it is
enough to study the distribution of $\lambda_{n,1}$. Deutsch,
Hildebrand and Wilf proved that the distribution of
$\sqrt{\frac{8}{n}}\ (\lambda_{n,1}-\frac{n}{2} )$ converges to the
distribution of the length of a random Gaussian vector in $\R^3$; in
other words $ \frac{8}{n} \left( \lambda_{n,1}-\frac{n}{2}
\right)^2 $ converges to the $\chi^2_3$ distribution with $3$
degrees of freedom (a careful reader may notice that the authors of
\cite{DeutschHildebrandWilf} use a non-standard definition of the
$\chi^2$ distributions and therefore they claim that
$\sqrt{\frac{8}{n}}\ (\lambda_{n,1}-\frac{n}{2} ) $ itself
converges to $\chi^2_3$). Their proof was based on an explicit
calculation of the number of the permutations which correspond to a
prescribed Young diagram with at most two rows.

\subsection{Case $d=\infty$}
\label{sub:Case-d=inf}
Another extreme of this problem is to
consider $d=\infty$; in other words, not to impose any restrictions
on the random permutations $\pi_n$. In this case the random Young
diagram $\lambda_n$ is distributed according to the Plancherel
measure. The authors of
\cite{BaikDeiftJohanssonDist1,Okounkov2000randompermutations,Johansson01}
proved that the joint distribution of longest rows of $\lambda$
(after appropriate rescaling) converges to the same distribution
(called Tracy--Widom distribution) as the joint distribution of the
biggest eigenvalues of a large random matrix from the Gaussian
Unitary Ensemble.


\subsection{The main result: intermediate values of $d$}

We equip the vector space of $d\times d$ Hermitian matrices with a
Gaussian probability measure with a density
$$ \frac{1}{Z_d}  e^{-\frac{1}{2} \Tr H^2} $$
with respect to the Lebesgue measure, where $Z_d$ is the normalizing constant. 
We say that a random matrix
$(A_{ij})_{1\leq i,j\leq d}$ distributed accordingly to this measure
is a Gaussian Unitary Ensemble (GUE) random matrix.


We call $B=A - \frac{1}{d} \Tr A$ a \emph{traceless Gaussian Unitary
Ensemble (GUE$_0$) random matrix}; it corresponds to the Gaussian
probability measure on the set of $d\times d$ Hermitian matrices
with trace zero and the density
$$ \frac{1}{Z'_d} e^{-\frac{1}{2} \Tr H^2} $$
with respect to the Lebesgue measure, where $Z'_d$ is the normalizing constant.


The joint distribution of eigenvalues for GUE is well-known
\cite{Mehta}, which allows us to find the corresponding distribution
for GUE$_0$; namely, if $x_1\geq\cdots\geq x_d$ are the eigenvalues of a
GUE$_0$ random matrix, then their joint distribution is supported on
the hyperplane $x_1+\cdots+x_d=0$ with the density
\begin{equation}
\label{eq:eigenvalues} 
\frac{1}{C_d}  e^{-\frac{x_1^2+\cdots+x_d^2}{2}} \prod_{i<j} (x_i-x_j)^2  
\end{equation}
with respect to the Lebesgue measure, where $C_d$ is the normalization
constant.

\begin{theorem}[Main theorem]
\label{theo:main} Let the integer $d\geq 1$ be fixed; for each $n\geq 1$ let
$\lambda_n=(\lambda_{n,1},\dots,\lambda_{n,d})$ be, as in Section
\ref{sub:Formula-problem}, a random Young diagram with $n$ boxes and with at
most $d$ rows.

Then the joint distribution of the random variables $ \left(
\sqrt{\frac{2d}{n}}\ (\lambda_{n,i}- \frac{n}{d}) \right)_{1\leq
i\leq d}$ converges, as $n\to\infty$, to the joint distribution of
the eigenvalues of a GUE$_0$ random matrix.
\end{theorem}
We postpone its proof to Section \ref{sec:Proof-main-result}.

\begin{corollary}
\label{cor:corollary}
Let $d\geq 1$ be fixed, and for each $n\geq 1$ let $\pi_n\in S_n$ be a random
permutation with the longest decreasing subsequence of length at most $d$. We
denote by $\lambda_{n,1}$ the length of its longest increasing subsequence. Then
the distribution of $ \sqrt{\frac{2d}{n}}\ (\lambda_{n,1}-\frac{n}{d}
) $ converges to the distribution of the largest eigenvalue of the
GUE$_0$ random matrix.
\end{corollary}

It should be pointed out that the distibution of eigenvalues of a GUE$_0$ random matrix appears also in a related  asymptotic problem \cite{Johansson01} of the distribution of the rows of a Young diagram associated (via RSK algorithm) to a random word consisting of $n$ letters in an alphabet of $d$ symbols in the limit of $n\to\infty$.

%

\subsection{Case $d=2$ revisited}

The set of $2\times 2$ Hermitian matrices with trace zero can be
viewed as a $3$-dimensional Euclidean space with a scalar product
$\langle A,B\rangle=\Tr AB$. A GUE$_0$ random matrix can
be viewed under this correspondence as a Gaussian random vector in
$\R^3$ the coordinates of which are independent with mean zero and 
variance $1$. Each $2\times 2$ Hermitian traceless matrix $A$ has
two eigenvalues $x_1=\lambda$, $x_2=-\lambda$, where
$\lambda=\frac{1}{\sqrt{2}} \|A\|=\sqrt{\frac{\langle A,A\rangle}{2}}$.
Therefore,
for a GUE$_0$ random matrix the corresponding random variable $2 (x_1)^2$ is
distributed like the square of the length of a standard Gaussian
random vector in $\R^3$, which coincides with the $\chi^2_3$
distribution; thus for $d=2$ Corollary \ref{cor:corollary} allows us to
recover the result of Deutsch,
Hildebrand and Wilf \cite{DeutschHildebrandWilf}.

\subsection{Idea of the proof}
In Section \ref{sec:Proof-main-result} we will prove Theorem \ref{theo:main},
the main result of this article.
Our proof will be based on an
explicit calculation of the number of standard Young tableaux with a prescribed
shape. The standard method to do this would be to use the hook-length formula,
which would be not convenient for our purposes. Instead, we will use the
determinantal formula of Frobenius and MacMahon. In order to make the
connection to random matrices more explicit we shall recall its proof due to
Zeilberger \cite{ZeilbergerAndreReflection} which is based on the observation
that a Young tableaux with at most $d$
rows can be viewed as a certain trajectory of $d$ non-colliding particles on a
line. Thus we will find explicitly the asymptotic joint distribution of the
rows of a Young diagram; this distribution turns out to coincide with
the distribution \eqref{eq:eigenvalues} of the eigenvalues of a GUE$_0$ random
matrix.

The reader may wonder if the connection between Young diagrams and random
matrices given by Theorem \ref{theo:main} might be purely accidental. In the
following paragraph we will argue why it is not the case and how deep
connections between Young diagrams and random matrices may be seen in our
proof of Theorem \ref{theo:main}.

%
%

In the above discussion we treated the distribution \eqref{eq:eigenvalues} of
the eigenvalues of a GUE$_0$ random matrix as granted; now let us think for a
moment about its derivation. GUE$_0$ is a Gaussian matrix; for this reason (up
to a simple scaling factor) it can be viewed as a value at some fixed time of a
matrix-valued Brownian bridge. It is known \cite{Dyson62,Grabiner99} that the
eigenvalues of a matrix-valued Brownion motion behave like Brownian motions
conditioned not to collide. Since a matrix-valued Brownian bridge is a
matrix-valued Brownian motion conditioned to be zero at time $1$, it follows
that its eigenvalues form Brownian motions conditioned not to collide and to be
zero at time $1$; in other words these eigenvalues form Brownian bridges
conditioned not to collide. In this way
the determinantal formula of Karlin and McGregor \cite{KarlinMcGregor59} can be
applied. In the conditioning procedure we assume that the original positions of
$d$ non-colliding particles are all different and we consider the limit as
these initial positions converge to zero; in this way their final distribution
is given by a continuous analogue of the formulas \eqref{eq:iloscdiagramow} and
\eqref{eq:calculation} which give the square of the number of Young tableaux of
a given shape, with the
transition probabilities replaced by the Gaussian kernels. One can easily check
that such a derivation of the  distribution of eigenvalues of a GUE$_0$ random
matrix follows \eqref{eq:calculation} very closely.

To summarize: our proof of the main result will be based on the
observation that both Young tableaux and the eigenvalues of matrix-valued
Brownian motions can be interpreted as non-colliding particles and applying
the determinantal formula of Karlin and McGregor \cite{KarlinMcGregor59}.

\subsection{Final remarks} We can see that both the case when $d$ is
finite and the case considered in Section \ref{sub:Case-d=inf}
corresponding to $d=\infty$ are asymptotically described by GUE
random matrices. It would be very interesting to find a direct link
between these two cases.


\section{Proof of the main result}
\label{sec:Proof-main-result}



%
%
%
%
%
%

For a function $f:\R\rightarrow\R$ we define its difference $\Delta_n
f:\R\rightarrow\R$ by
$$ \Delta_n f(y)= \frac{ f\left(y+\sqrt{\frac{d}{n}}  \right) -f(y)
}{\sqrt{\frac{d}{n}}}. $$
By iterating we define $\Delta_n^\alpha f$ for any integer
$\alpha\geq 0$.
We also define its shift $S_n f:\R\rightarrow \R$ by
$$ S_n f(y)= f\left(y+\sqrt{\frac{d}{n}}\right). $$
Notice that $S_n^{\alpha}f$ is well-defined for any integer $\alpha$.

\begin{lemma}
\label{lem:PoissonOk}
For each $n$ we define a function $f_n:\R\rightarrow\R$ which is constant on
each
interval of the form $\left[\frac{k- \frac{n}{d} }{\sqrt{\frac{n}{d}} }
, \frac{k+1- \frac{n}{d} }{\sqrt{\frac{n}{d}} }\right)$ for each integer $k$
and such that
\begin{equation} 
\label{eq:Poisson} f_n\left( \frac{k- \frac{n}{d} }{\sqrt{\frac{n}{d}} }
\right)= \begin{cases} \sqrt{\frac{n}{d}} \frac{\left(\frac{n}{d}\right)^k
e^{-\frac{n}{d}} }{k!} & \text{if $k$ is a non-negative integer,} \\
0 & \text{if $k$ is a negative integer.}
\end{cases}
\end{equation}

Then for each integer $\alpha\geq 0$ and $y\in\R$
\begin{equation} 
\label{eq:zbieznosc}
\lim_{n\to\infty }  \Delta^\alpha_n
f_n (y)  = \frac{d^\alpha}{dy^\alpha} \frac{1}{\sqrt{2\pi}}
e^{\frac{-y^2}{2}}. 
\end{equation}

Furthermore, for each $\alpha\geq 0$ there exists a polynomial $P_{\alpha}$ such
that \begin{equation} 
\label{eq:dominacja}
\big| \Delta^\alpha_n
f_n (y) \big| < P_{\alpha}(y) e^{-|y|} 
\end{equation}
holds true for all $n$ and $y$.

\end{lemma}
\begin{proof}
Before presenting the proof we notice that $f_n$ is a density of a probability
measure arising as follows: we normalize the Poisson distribution  with
the parameter $\nu=\frac{n}{d}$ in order to have mean $0$ and variance $1$ and
we convolve it with a uniform distribution on the interval $\left[ 0,
\sqrt{\frac{d}{n}} \right] $;
therefore
\eqref{eq:zbieznosc} states for $\alpha=0$ that the suitably rescaled
probabilities of the Poisson distribution converge to the density of the normal
distribution. The case $\alpha\geq 1$ shows that this convergence holds true
also for differences (respectively, derivatives). 

The proof of \eqref{eq:zbieznosc} in the case $\alpha=0$ is a straightforward
application of the Stirling approximation $\log z!=\left(z+\frac{1}{2} \right) \log z - z +\frac{\log 2\pi}{2} + O(z^{-1})$, namely for $y=\frac{k- \frac{n}{d} }{\sqrt{\frac{n}{d}} }$ such that $k$ is an integer we denote $c=\frac{n}{d}$. Then
\begin{multline*} \log f_n(y)= \left( c+ y \sqrt{c} + \frac{1}{2} \right) \log c- c - \log \left( c + y \sqrt{c} \right)! = \\  
-\left( c + y \sqrt{c} +\frac{1}{2} \right) \log \left( 1 + \frac{y}{ \sqrt{c}} \right) + y \sqrt{c} - \frac{\log 2\pi}{2} + O\left({c}^{-1}\right)
 =\\ -\frac{y^2}{2} -\frac{\log 2\pi}{2} + O\left({c}^{-\frac{1}{2}}\right),
\end{multline*}
where the above equalities hold true asymptotically for $y$ bounded and $c\to\infty$.

In order to treat the case $\alpha\geq 1$ we observe
that the iterated derivative on the right-hand side of \eqref{eq:zbieznosc} can
be calculated by using the following three rules:
$$ \frac{d}{dy} e^{\frac{-y^2}{2}}= -y e^{\frac{-y^2}{2}}; \qquad \frac{d}{dy} y
= 1; \qquad \frac{d}{dy} (\phi \psi)= \left(\frac{d}{dy}\phi \right) \psi + \phi
\frac{d}{dy}\psi. $$
Similarly, the iterated difference on the left-hand side of
\eqref{eq:zbieznosc} can be calculated using the following three rules:
$$ \Delta_n f_n = - g_n S f_n; \qquad 
\Delta_n g_n=1; \qquad \Delta_n (ab)= (\Delta_n a) b+ (S_n a) \Delta_n b   $$
where $g_n:\R\rightarrow\R$ is a function which is constant on each
interval of the form $\left[\frac{k- \frac{n}{d} }{\sqrt{\frac{n}{d}} }
, \frac{k+1- \frac{n}{d} }{\sqrt{\frac{n}{d}} }\right)$ for each integer $k$
and such that
$$ g_n\left( \frac{k- \frac{n}{d} }{\sqrt{\frac{n}{d}} } \right)=
\frac{k+1-\frac{n}{d}}{\sqrt{\frac{n}{d}}}. $$

For each integer $\beta$ we have $\lim_{n\to\infty} (S_n^\beta f_n)(y) =
\frac{1}{\sqrt{2\pi}} e^{-\frac{y^2}{2}} $ and $\lim_{n\to\infty}
(S_n^{\beta}g_n)(y)=y $ therefore each term contributing to the left-hand
side of \eqref{eq:zbieznosc} converges to its counterpart on the right-hand
side of \eqref{eq:zbieznosc}, which finishes the proof of \eqref{eq:zbieznosc}.

We consider $y=\frac{k- \frac{n}{d} }{\sqrt{\frac{n}{d}} }$; then
\begin{equation}
\label{eq:logarytmy}
\frac{\log f_n\left( \frac{k-\frac{n}{d}}{\sqrt{\frac{n}{d}}}  \right)- 
\log f_n\left( \frac{k-1-\frac{n}{d}}{\sqrt{\frac{n}{d}}}  \right)
}{\sqrt{\frac{n}{d}}} = - \log \left( 1+ \frac{y}{\sqrt{\frac{n}{d}}} \right) \sqrt{\frac{n}{d}}. 
\end{equation}
There is a constant $C_1<0$ with a property that if $y<C_1$ then the right-hand side of \eqref{eq:logarytmy} is greater than $1$ for any value of $n$. It follows that if $y_i=\frac{k_i- \frac{n}{d} }{\sqrt{\frac{n}{d}} }$ for $i\in\{1,2\}$ and $y_1<y_2\leq C_1$ then
\begin{equation}
\label{eq:szacowanielogarytmu1}
 f_n(y_1)\leq f_n(y_2) e^{y_1-y_2}. 
\end{equation}
Similarly we find a constant $C_2>0$ with a property that if $C_2\leq y_1<y_2$ then
\begin{equation}
\label{eq:szacowanielogarytmu2}
f_n(y_2)\leq f_n(y_1) e^{y_1-y_2}. 
\end{equation}
For $\alpha=0$ inequality \eqref{eq:dominacja} holds true for $y$ in a small neighborhood of the interval $[C_1,C_2]$ for $P_{\alpha}$ being a sufficiently big constant which follows from \eqref{eq:zbieznosc} and compactness argument. Inequality \eqref{eq:dominacja} holds true outside of the interval $[C_1,C_2]$ by inequalities \eqref{eq:szacowanielogarytmu1} and \eqref{eq:szacowanielogarytmu2}.

The case $\alpha\geq 1$ can be proved in an analogous way to the above proof of
\eqref{eq:zbieznosc}: we show that $\Delta^\alpha_n f_n$ is a sum of the terms
of the form $(S_n^{\beta_1} g_n) \cdots (S_n^{\beta_l} g_n)(S^{\beta}_n f_n)$
and the absolute value of each such a term can be easily bounded by $P(y)
e^{-|y|}$, where $P$ is a suitably chosen polynomial. 


\end{proof}

\begin{proof}[Proof of Theorem \ref{theo:main}]
The following discussion is based on the work of Zeilberger
\cite{ZeilbergerAndreReflection}. Every Young tableau $T$ with at most $d$ rows
and $n$ boxes can be interpreted
as a trajectory of $d$ non-colliding particles $x_1(t),\dots,x_d(t)$ on the real
line as follows.
We set 
$$x_i(t)= d+1-i+(\text{number of boxes of $T$ in row $i$ which are not bigger
than $t$)}. $$
In other words: the initial positions of the particles are given by
$\big(x_1(0),\dots,x_d(0)\big)=(d,d-1,\dots,1)$. In each step one of the
particles jumps to the right; the number of the particle which jumps in step
$t$ is equal to the number of the row of the Young diagram $T$ which carries the
box with a label $t$. The condition that $T$ is a standard Young tableau is
equivalent to $x_1(t)>\cdots>x_d(t)$ for every value of $0\leq t\leq n$. 

Thus
the results of Karlin and McGregor \cite{KarlinMcGregor59} can be applied and
the number of standard Young tableaux of the shape $\lambda_1,\dots,\lambda_d$,
where $|\lambda|=\lambda_1+\cdots+\lambda_d=n$, is equal to the determinant
\begin{multline}
\label{eq:calculation} 
N_{\lambda_1,\dots,\lambda_n}= n!
\begin{vmatrix} 
\frac{1}{\lambda_1!} & \frac{1}{(\lambda_1+1)!}  & \cdots &  \frac{1}{(\lambda_1+d-1)!} \\
\frac{1}{(\lambda_2-1)!} & \frac{1}{\lambda_2!} & \cdots & \frac{1}{(\lambda_2+d-2)!} \\
\vdots & \vdots & \ddots & \vdots \\ \frac{1}{(\lambda_d-d+1)!} &
\frac{1}{(\lambda_d-d+2)!}  & \cdots &\frac{1}{\lambda_d!}  \end{vmatrix} = \\
\frac{n! e^n}{\left(\frac{n}{d} \right)^{n+\frac{d}{2}} } \begin{vmatrix} 
f_{n}(y_1) & S_n f_n(y_1) & \cdots & S^{d-1 }_n f_n(y_1) \\
S_n^{-1} f_{n}(y_2) & f_n(y_2) & \cdots & S_n^{d-2} f_n(y_2) \\
\vdots & \vdots & \ddots & \vdots \\
S_n^{-d+1} f_{n}(y_d) & S^{-d+2}_n f_n(y_d) & \cdots & f_n(y_d)
\end{vmatrix} = \\
\frac{n! e^n}{\left(\frac{n}{d} \right)^{n+\frac{d (d+1)}{4}} } 
\begin{vmatrix} 
f_{n}(y_1) &  \Delta_n f_n(y_1)             
& \cdots & \Delta^{d-1}_n f_n(y_1) \\
S_n^{-1}f_{n} (y_2) & \Delta_n S_n^{-1} f_n(y_2) & \cdots &
\Delta_n^{d-1} S_n^{-1} f_n(y_2) \\
\vdots & \vdots & \ddots & \vdots \\
S_n^{-d+1} f_{n}(y_d) &  \Delta_n S_n^{-d+1} f_n(y_d) & \cdots &
\Delta_n^{d-1} S_n^{-d+1} f_n(y_d)
\end{vmatrix},
\end{multline} 
where 
\begin{equation} 
\label{eq:y}
y_i= \frac{\lambda_i- \frac{n}{d} }{\sqrt{\frac{n}{d}}}.
\end{equation}

We are interested in a probability distribution on Young diagrams with $n$
boxes with the probability of $(\lambda_1,\dots,\lambda_d)$ equal to
\begin{equation}
\label{eq:iloscdiagramow}
 \frac{1}{C_{n,d}} (N_{\lambda_1,\dots,\lambda_d})^2, 
\end{equation}
where $C_{n,d}$ is the suitably chosen normalizing constant. 
Clearly,
\begin{multline*} C_{n,d}  \frac{\left(\frac{n}{d}
\right)^{2n+\frac{d^2+2d-1}{2} }}{(n!)^2 e^{2n}}  = \\
{\sum_{\lambda_1,\dots,\lambda_{d-1}} \left( \sqrt{\frac{n}{d}}
\right)^{d-1} }  
\begin{vmatrix} 
f_{n}(y_1) &  \Delta_n f_n(y_1)            
& \cdots & \Delta^{d-1}_n f_n(y_1) \\
S_n^{-1}f_{n} (y_2) & \Delta_n S_n^{-1} f_n(y_2) & \cdots &
\Delta_n^{d-1} S_n^{-1} f_n(y_2) \\
\vdots & \vdots & \ddots & \vdots \\
S_n^{-d+1} f_{n}(y_d) &  \Delta_n S_n^{-d+1} f_n(y_d) & \cdots &
\Delta_n^{d-1} S_n^{-d+1} f_n(y_d)
\end{vmatrix}^2,
\end{multline*}
where the sum runs over $\lambda_1,\dots,\lambda_{d-1}$  such that
for $\lambda_{d}=n-(\lambda_1+\cdots+\lambda_{d-1})$ we have that
$\lambda_1,\dots,\lambda_d$ is a Young diagram with $n$ boxes.
The right-hand side can be viewed as a Riemann sum; Lemma
\ref{lem:PoissonOk} shows that 
the dominated convergence theorem can be applied 
(with the dominating function of the form $P(y_1,\dots,y_d) e^{-2 (|y_1|+\cdots+|y_d|)}$, where $P$ is some polynomial)
and 
\begin{multline*} \lim_{n\to\infty} C_{n,d}  \frac{\left(\frac{n}{d}
\right)^{2n+\frac{d^2+2d-1}{2} }}{(n!)^2 e^{2n}}  = \\ \int_{y_1,\dots,y_{d-1}}
\begin{vmatrix}
e^{-\frac{y_1^2}{2}} & \frac{d}{dy_1} e^{-\frac{y_1^2}{2}} & \cdots &
\frac{d^{d-1}}{dy_1^{d-1}} e^{-\frac{y_1^2}{2}} \\ 
e^{-\frac{y_2^2}{2}} & \frac{d}{dy_2} e^{-\frac{y_2^2}{2}} & \cdots &
\frac{d^{d-1}}{dy_2^{d-1}} e^{-\frac{y_2^2}{2}} \\
\vdots & \vdots & \ddots & \vdots \\
e^{-\frac{y_d^2}{2}} & \frac{d}{dy_d} e^{-\frac{y_d^2}{2}} & \cdots &
\frac{d^{d-1}}{dy_d^{d-1}} e^{-\frac{y_d^2}{2}}  
\end{vmatrix}^2 dy_1 \cdots dy_{d-1},
\end{multline*}
where the integral runs over $(y_1,\dots,y_{d-1})$ such that for
$y_d=-(y_1+\cdots+y_{d-1})$ we have $y_1\geq \cdots \geq y_d$.

Since the limit density defines a probability measure, in the
limit $n\to\infty$ the random variables $(y_1,\dots,y_{d-1})$
(please notice that due to the constraint $y_1+\cdots+y_d=0$ the value of $y_d$
is uniquely determined by $y_1,\dots,y_{d-1}$) converge in distribution to the
probability measure on the set $y_1\geq y_2\geq \cdots \geq y_{d-1} \geq
-(y_1+\cdots+y_{d-1})$ with a density
\begin{multline*}\frac{1}{C'_d} \begin{vmatrix}
e^{-\frac{y_1^2}{2}} & \frac{d}{dy_1} e^{-\frac{y_1^2}{2}} & \cdots &
\frac{d^{d-1}}{dy_1^{d-1}} e^{-\frac{y_1^2}{2}} \\ 
e^{-\frac{y_2^2}{2}} & \frac{d}{dy_2} e^{-\frac{y_2^2}{2}} & \cdots &
\frac{d^{d-1}}{dy_2^{d-1}} e^{-\frac{y_2^2}{2}} \\
\vdots & \vdots & \ddots & \vdots \\
e^{-\frac{y_d^2}{2}} & \frac{d}{dy_d} e^{-\frac{y_d^2}{2}} & \cdots &
\frac{d^{d-1}}{dy_d^{d-1}} e^{-\frac{y_d^2}{2}}  
\end{vmatrix}^2  = \\
\frac{1}{C'_d} \begin{vmatrix}
p_0(y_1) e^{-\frac{y_1^2}{2}} & p_1(y_1) e^{-\frac{y_1^2}{2}} & \cdots &
p_{d-1}(y_1) e^{-\frac{y_1^2}{2}} \\ 
p_0(y_2) e^{-\frac{y_2^2}{2}} & p_1(y_2) e^{-\frac{y_2^2}{2}} & \cdots &
p_{d-1}(y_2) e^{-\frac{y_2^2}{2}} \\
\vdots & \vdots & \ddots & \vdots \\
p_0(y_d) e^{-\frac{y_d^2}{2}} & p_1(y_d) e^{-\frac{y_d^2}{2}} & \cdots &
p_{d-1}(y_d) e^{-\frac{y_d^2}{2}}  
\end{vmatrix}^2  \end{multline*}
for a suitably chosen normalizing constant $C'_d$, where $\frac{d^k}{dz^{k}} e^{-\frac{z^2}{2}}=p_k(z)   e^{-\frac{z^2}{2}}$ for some
polynomial $p_k$ (related to Hermite polynomials). Since $p_k(z)=(-z)^k + \text{(summands of lower degree)}$ the above expression takes a simpler form:
\begin{multline*}\frac{1}{C'_d} \begin{vmatrix}
 e^{-\frac{y_1^2}{2}} & (-y_1) e^{-\frac{y_1^2}{2}} & \cdots &
(-y_1)^{d-1} e^{-\frac{y_1^2}{2}} \\ 
 e^{-\frac{y_2^2}{2}} & (-y_2) e^{-\frac{y_2^2}{2}} & \cdots &
(-y_2)^{d-1} e^{-\frac{y_2^2}{2}} \\
\vdots & \vdots & \ddots & \vdots \\
 e^{-\frac{y_d^2}{2}} & (-y_d) e^{-\frac{y_d^2}{2}} & \cdots &
(-y_d)^{d-1} e^{-\frac{y_d^2}{2}}  
\end{vmatrix}^2  = \\
\frac{1}{C'_d} e^{-(y_1^2+\cdots+y_d^2)  } \prod_{1\leq
i<j\leq d} (y_i-y_j)^2. 
\end{multline*}

 When we set $x_i=\sqrt{2}
y_i= \sqrt{\frac{2d}{n}}\ (\lambda_{n,i}- \frac{n}{d})$ it becomes clear that
the limit distribution of $(x_1,\dots,x_d)$ coincides with the distribution
\eqref{eq:eigenvalues} of 
the eigenvalues of a $GUE_0$ random matrix, which finishes the  proof.
\end{proof}

\section{Acknowledgements}
\label{sec:Acknowl} 
The research was performed during a visit to Queens University. I thank Jonathan
Novak for many discussions and for pointing out the reference
\cite{DeutschHildebrandWilf}. I thank Roland Speicher for invitation and
hospitality during this stay.

Research supported by the MNiSW research grant 1 P03A 013 30, by the EU Research
Training Network `QP-Applications', contract HPRN-CT-2002-00279 and  by the EC
Marie Curie Host Fellowship for the Transfer of Knowledge `Harmonic Analysis,
Nonlinear Analysis and Probability', contract MTKD-CT-2004-013389.

 \bibliographystyle{alpha}
 \bibliography{biblio}

\begin{thebibliography}{DHW03}

\bibitem[BDJ99]{BaikDeiftJohanssonDist1}
Jinho Baik, Percy Deift, and Kurt Johansson.
\newblock On the distribution of the length of the longest increasing
  subsequence of random permutations.
\newblock {\em J. Amer. Math. Soc.}, 12(4):1119--1178, 1999.

\bibitem[DHW03]{DeutschHildebrandWilf}
Emeric Deutsch, A.~J. Hildebrand, and Herbert~S. Wilf.
\newblock Longest increasing subsequences in pattern-restricted permutations.
\newblock {\em Electron. J. Combin.}, 9(2):Research paper 12, 8 pp.
  (electronic), 2002/03.
\newblock Permutation patterns (Otago, 2003).

\bibitem[Dys62]{Dyson62}
Freeman~J. Dyson.
\newblock A {B}rownian-motion model for the eigenvalues of a random matrix.
\newblock {\em J. Mathematical Phys.}, 3:1191--1198, 1962.

\bibitem[Gra99]{Grabiner99}
David~J. Grabiner.
\newblock Brownian motion in a {W}eyl chamber, non-colliding particles, and
  random matrices.
\newblock {\em Ann. Inst. H. Poincar\'e Probab. Statist.}, 35(2):177--204,
  1999.

\bibitem[Joh01]{Johansson01}
Kurt Johansson.
\newblock Discrete orthogonal polynomial ensembles and the {P}lancherel
  measure.
\newblock {\em Ann. of Math. (2)}, 153(1):259--296, 2001.

\bibitem[KM59]{KarlinMcGregor59}
Samuel Karlin and James McGregor.
\newblock Coincidence probabilities.
\newblock {\em Pacific J. Math.}, 9:1141--1164, 1959.

\bibitem[Meh91]{Mehta}
Madan~Lal Mehta.
\newblock {\em Random matrices}.
\newblock Academic Press Inc., Boston, MA, second edition, 1991.

\bibitem[Oko00]{Okounkov2000randompermutations}
Andrei Okounkov.
\newblock Random matrices and random permutations.
\newblock {\em Internat. Math. Res. Notices}, (20):1043--1095, 2000.

\bibitem[Zei83]{ZeilbergerAndreReflection}
Doron Zeilberger.
\newblock Andr\'e's reflection proof generalized to the many-candidate ballot
  problem.
\newblock {\em Discrete Math.}, 44(3):325--326, 1983.

\end{thebibliography}

\end{document}